\documentclass[12pt, reqno]{amsart}

\usepackage{amsfonts,amssymb,latexsym,amsmath, amsxtra}
\usepackage{verbatim}

\usepackage{amsthm,amssymb,amstext,amscd,amsfonts,amsbsy,amsxtra,latexsym,amsmath}
\usepackage[english]{babel}
\usepackage[latin1]{inputenc}
\usepackage{graphicx}
\usepackage[all]{xy}
\numberwithin{figure}{section}

\newcommand{\field}[1]{\mathbb{#1}}
\newcommand{\N}{\field{N}}
\newcommand{\Z}{\field{Z}}

\newcommand{\C}{\field{C}}

\renewcommand{\S}{\mathcal{S}}
\newcommand{\GG}{\mathcal{GG}}
\newcommand{\G}{\mathcal{G}}
\newcommand{\calT}{\mathcal{T}}
\newcommand{\qbinom}[2]{\genfrac{[}{]}{0pt}{}{#1}{#2}}

\theoremstyle{plain}
\newtheorem{theorem}{Theorem}[section]
\newtheorem*{theorem*}{Theorem}

\newtheorem{lemma}[theorem]{Lemma}

\newtheorem*{conjecture*}{Conjecture}

\theoremstyle{definition}

\theoremstyle{remark}
\newtheorem{remark}{Remark}

\newtheorem*{remark*}{Remark}
\newtheorem*{remarks*}{Remarks}

\numberwithin{equation}{section}

\begin{document}

\title[Double series for Schur partitions]{Double series representations for Schur's partition function and related identities}
\author{George Andrews}
\address{Department of Mathematics\\The Pennsylvania State University \\ University Park, PA 16802 \\ U.S.A.}
\email{gea1@psu.edu}
\author{Kathrin Bringmann}
\address{Mathematical Institute\\University of
Cologne\\ Weyertal 86-90 \\ 50931 Cologne \\Germany}
\email{kbringma@math.uni-koeln.de}
\author{Karl Mahlburg}
\address{Department of Mathematics \\
Louisiana State University \\
Baton Rouge, LA 70802\\ U.S.A.}
\email{mahlburg@math.lsu.edu}


\date{\today}


\thanks{The first author was partially supported by NSA grant H98230-12-1-0205. The research of the second author was supported by the Alfried Krupp Prize for Young University Teachers of the Krupp Foundation and the research leading to these results has received funding from the European Research Council under the European Union's Seventh Framework Programme (FP/2007-2013) / ERC Grant agreement n. 335220 - AQSER.  The third author was supported by NSF Grant DMS-1201435.}

\begin{abstract}
We prove new double summation hypergeometric $q$-series representations for several families of partitions, including those that appear in the famous product identities of G\"ollnitz, Gordon, and Schur. We give several different proofs for our results, using bijective partitions mappings and modular diagrams, the theory of $q$-difference equations and recurrences, and the theories of summation and transformation for $q$-series. We also consider a general family of similar double series and highlight a number of other interesting special cases.
\end{abstract}

\maketitle

\section{Introduction and statement of results}

For an integer partition $\lambda$ with parts $\lambda_1 \geq \cdots \geq \lambda_\ell$, we denote its {\it size} by $|\lambda| := \lambda_1 + \dots + \lambda_\ell$, and its {\it length} (or number of parts) by $\ell(\lambda)$ (see \cite{And98} for the basic theory of partitions).
The first family that we consider is the {\it Schur partitions}. Let $\S$ denote the set of partitions $\lambda$ such that $\lambda_i - \lambda_{i+1} \geq 3$, with strict inequality if either part is a multiple of $3$.
We define the two-parameter generating function for Schur's partitions as
\begin{equation*}
f_{\S}(x) = f_{\S}(x;q) := \sum_{\lambda \in \S} x^{\ell(\lambda)} q^{|\lambda|}.
\end{equation*}
Note that this function was previously denoted by $f_0(x)$ in \cite{And68} and as $f_{3,1}(x;q)$ in \cite{BLM13}.
It is of interest due to a striking infinite product identity, as Schur's Second Partition Theorem \cite{Sch26} states that
\begin{equation}
\label{E:Schur=prod}
f_{\S}(1;q) = \left(-q; q^3\right)_\infty \left(-q^2; q^3\right)_\infty.
\end{equation}
Here and throughout the paper, we adopt standard notation for hypergeometric $q$-series, which is reviewed in detail in Section \ref{S:qseries}.

Our first result gives a new double series representation for the two-parameter generating function for Schur partitions.
\begin{theorem}
\label{T:fSdouble}
We have
\begin{equation*}
f_{\S}(x) = \sum_{m, n \geq 0} \frac{(-1)^n x^{m + 2n} q^{3n(3n+2m) + \frac{m(3m-1)}{2}}}{(q;q)_m (q^6; q^6)_n}.
\end{equation*}
\end{theorem}
\begin{remark*}
Theorem \ref{T:fSdouble} is particularly notable because, unlike some of the other examples that are discussed in this paper, there are not any simple combinatorial $q$-series representations for $f_{\S}(x)$ that appear in the literature. For example, one of the simplest expressions is (2.15) of \cite{And68}, which states that
\begin{equation*}
f_{\S}(x) = \left(x; q^3\right)_\infty \sum_{n \geq 0}
\frac{x^n \left(-q, -q^2; q^3\right)_n}{\left(q^3; q^3\right)_n}.
\end{equation*}
One significant disadvantage of this expression is the fact that it is an indeterminate limit when $x = 1$, whereas our new double series does not display such a singularity.
\end{remark*}

\begin{remark*}
Sills proved a similar double series representation for the product that appears in the notable partition identity of Capparelli, which was originally conjectured in \cite{Cap88}. Equation (1.4) of \cite{Sil04} states that
\begin{equation*}
\sum_{n \geq 0} \sum_{j = 0}^{2n} \frac{\left(\frac{n-j+1}{3}\right)q^{n^2}}{(q;q)_{2n-j} (q;q)_j}
= \frac{1}{\left(q^2, q^3, q^9, q^{10}; q^{12}\right)_\infty},
\end{equation*}
where $(\frac{\bullet}{p})$ denotes the Legendre symbol.
As with Schur's identity \eqref{E:Schur=prod}, the above product is equivalent to the generating function for partitions whose successive parts satisfy certain gap conditions.
Sills' proof uses the theory of Bailey pairs, which does not appear in any of our proofs.
\end{remark*}

We next consider partitions studied by G\"ollnitz and Gordon.
Let $\GG$ denote the set of {\it G\"ollnitz-Gordon partitions}, which are those partitions satisfying $\lambda_i - \lambda_{i+1} \geq 2$, with strict inequality if either part is even. Furthermore, let $\GG_t$ denote those partitions in $\GG$ with all parts at least $t$. A direct combinatorial argument (conditioning on the number of parts) shows that
\begin{equation}
\label{E:Gsum}
f_{\GG}(x;q) := \sum_{\lambda \in \GG} x^{\ell(\lambda)} q^{|\lambda|}
= \sum_{n \geq 0} \frac{x^n q^{n^2} (-q; q^2)_n}{(q^2; q^2)_n}. \\
\end{equation}
Furthermore, if the sum is instead taken over partitions $\lambda \in \GG_t$ and $t$ is odd, it is clear that the resulting generating function is $f_{\GG}(xq^{t-1};q)$.
The G\"ollnitz-Gordon identities, which were independently proven in \cite{Gol67,Gor65}, then state that for $t = 1$ or $3$, we have the following product formulas:
\begin{equation*}
\label{E:GGprod}
f_{\GG_t}\left(1;q\right) = \frac{1}{\left(q^t, q^4, q^{8-t}; q^8\right)_\infty}.
\end{equation*}

Similarly, let $\G$ denote the set of {\it G\"ollnitz partitions}, which are those partitions satisfying $\lambda_i - \lambda_{i+1} \geq 2$, where now the inequality is strict if either part is odd. Furthermore, let $\G_t$ denote the subset of such partitions where the smallest part is at least $t$. Denote the corresponding generating function by
\begin{equation*}
f_{\G}(x;q) := \sum_{\lambda \in \G} x^{\ell(\lambda)} q^{|\lambda|},
\end{equation*}
and let $f_{\G_t}$ be the sum over $\G_t$. As above, if $t$ is odd, then $f_{\G_t}(x;q) = f_{\G}(xq^{t-1};q)$.
G\"ollnitz specifically studied the cases $t=1, 2$, and noted the combinatorial formulas
\begin{align*}
f_{\G_1}(x;q) &= \sum_{n \geq 0} \frac{x^n q^{n^2+n} \left(-q^{-1}; q^2\right)_n}{\left(q^2; q^2\right)_n}, \\
f_{\G_2}(x;q) &= f_{\GG}(xq; q).
\end{align*}
The Little G\"ollnitz Theorem \cite{Gol67} gives the product identities
\begin{align}
\label{E:LG1}
f_{\G_1}(1;q) &= \frac{1}{\left(q, q^5, q^6; q^8\right)_\infty}, \\
\label{E:LG2}
f_{\G_2}(1;q) &= \frac{1}{\left(q^2, q^3, q^7; q^8\right)_\infty}.
\end{align}

Our next result gives new double series representations of the generating functions for G\"ollnitz-Gordon and G\"ollnitz partitions.
\begin{theorem}
\label{T:GG}
We have
\begin{align}
\label{E:GGdouble}
f_{\GG}(x;q) &= \sum_{k,m \geq 0} \frac{(-1)^k x^{m+2k} q^{m^2 + 4mk + 6k^2}}{(q;q)_m (q^4; q^4)_k}, \\
\label{E:Gdouble}
f_{\G}(x;q) &= \sum_{k,m \geq 0} \frac{(-1)^k x^{m+2k} q^{m^2 + 4mk + 6k^2 - 2k}}{(q;q)_m (q^4; q^4)_k}.
\end{align}
\end{theorem}
\begin{remark*}
Alladi and Berkovich also obtained alternative double series representations for series related to G\"ollnitz partitions. In particular, they proved (see (3.6) and (3.9) of \cite{AB05}) the identities
\begin{align}
\label{E:AB1}
\sum_{k,m \geq 0} \frac{w^k q^{m^2 + 2mk + 2k^2 + m-k}}{\left(q^2; q^2\right)_m \left(q^2; q^2\right)_k}
& = \left(-q^2, -wq^3, -q^4; q^4\right)_\infty, \\
\label{E:AB2}
\sum_{k,m \geq 0} \frac{w^k q^{m^2 + 2mk + 2k^2 + m+k}}{\left(q^2; q^2\right)_m \left(q^2; q^2\right)_k}
& = \left(-wq, -q^2, -q^4; q^4\right)_\infty.
\end{align}
If we set $w = 1$, these specialize to the products in \eqref{E:LG1} and \eqref{E:LG2}, respectively. Note that unlike our parameter $x$, which arises from partitions satisfying gap conditions (i.e., the ``sum side''), their parameter $w$ instead comes from partitions with congruential restrictions (the ``product side''). 

We also note that Alladi and Berkovich provided both combinatorial and analytic proofs of \eqref{E:AB1} and \eqref{E:AB2} in \cite{AB05}.
\end{remark*}

In order to further explore the role of double series in combinatorial partition identities, we define the general family
\begin{equation*}
\label{E:defineR}
R(s,t,\ell,u,v,w):= \sum_{n\geq 0}q^{s \binom{n} 2 + tn}r(\ell,u,v,w;n),
\end{equation*}
where
\begin{equation}
\label{E:definer}
r(\ell,u,v,w;n) := \sum_{j\geq0} \frac{(-1)^j q^{uv \frac{j(j-1)}{2}+(w-u\ell)j}}{\left(q;q\right)_{n-uj}(q^{uv};q^{uv})_{j}}.
\end{equation}
To illustrate this notation, we note that the double series associated to the partition identities described above may be written as
\begin{align}
\notag
R(3,t,0,2,3,4) &= f_\S\left(q^{t-1};q\right), \\
\label{E:R=GG}
R(2,t,0,2,2,2) &= f_\GG\left(q^{t-1};q\right), \\
\label{E:R=G}
R(2,t,1,2,2,2) &= f_{\G}\left(q^{t-1};q\right).
\end{align}

Our final result identifies several additional cases in which the double series is equivalent to an infinite product.
\begin{theorem}
\label{T:Rindentities}
 We have the following identities
\begin{align}
 R(2,1,1,1,2,2)&= \left(-q;q^2\right)_{\infty} \label{E:R=Euler1}, \\
 R(2,2,1,1,2,2)&= \left(-q^2;q^2\right)_{\infty} \label{E:R=Euler2}, \\
 R(m,m,1,1,1,2)&= \frac{\left(q^{2m};q^{2m}\right)_{\infty}}{\left(q^m;q^{2m}\right)_{\infty}}. \label{E:Rm=prod}
\end{align}
\end{theorem}
\begin{remark*}
Note that the double-series notation is not uniquely determined, as $$r(\ell,u,v,w;n) = r(0,u,v,w-u\ell;n).$$ However, we write the parameters in this way as it is often convenient to isolate $\ell$ (cf. Lemma \ref{L:recurrences}).
\end{remark*}

We give multiple proofs of our new double series identities, utilizing three primary techniques: combinatorial/bijective mappings of partitions, $q$-difference equations and analytic recurrences, and finally, the theory of hypergeometric $q$-series. Indeed, one of the goals of this paper is to show how these techniques can be adapted to solve similar problems.

The remainder of the paper is structured as follows. In the next section we review the basic notation of hypergeometric $q$-series, and collect a number of useful identities. In Section \ref{S:Comb} we give combinatorial proofs of Theorems \ref{T:fSdouble} and \ref{T:GG}. We turn to analytic methods in Section \ref{S:AnalyticI}, using $q$-difference equations to prove Theorems \ref{T:fSdouble} -- \ref{T:Rindentities}. Section \ref{S:AnalyticII} gives analytic proofs of Theorems \ref{T:GG} and \ref{T:Rindentities} using basic identities from the theory of hypergeometric $q$-series.

\section{Basic hypergeometric $q$-series}
\label{S:qseries}

In this section we recall the basic definitions and notation for hypergeometric $q$-series. We also record a number of identities that are useful in the evaluation of the generating functions that are the main topic of the paper.

If $a \in \C$ and $n \in \Z$, the {\it $q$-Pochhammer} symbol is defined by setting
\begin{align*}
(a;q)_{\infty} & := \prod_{k \geq 0} \left(1 - aq^k\right), \\
(a;q)_n & := \frac{(a;q)_\infty}{\left(aq^n; q\right)_\infty}.
\end{align*}
Note that this definition gives the convenient convention that
\begin{equation}
\label{E:q_n<0}
\frac{1}{(q;q)_n} = 0 \qquad \text{if } n < 0.
\end{equation}
For example, this means that summations such as \eqref{E:definer} can be extended to all integers $j$.
We also adopt the shorthand notations $(a)_n := (a;q)_n$
 and $(a_1, \dots, a_r)_n := (a_1)_n \cdots (a_r)_n$.

Next, we note the identity
\begin{equation}
\label{E:qnegn}
\left(q^{-n}\right)_j = (-1)^j q^{-nj + \frac{j(j-1)}{2}} \frac{(q)_{n}}{(q)_{n-j}}.
\end{equation}
We also have the useful limit evaluation
\begin{equation}
\label{E:qtlim}
\lim_{t \rightarrow 0} t^j \left(-\frac{q}{t}\right)_j = q^{\frac{j(j+1)}{2}}.
\end{equation}

If $0 \leq m \leq n$, then the {\it $q$-binomial coefficient} is denoted by
\begin{equation*}
\qbinom{n}{m}_q := \frac{(q;q)_n}{(q;q)_m (q;q)_{n-m}}.
\end{equation*}
The $q$-Binomial Theorem (see Theorem 3.3 of \cite{And98}) states that if $n \geq 0$,
\begin{equation}
\label{E:qbinom}
(a)_n = \sum_{j = 0}^n \qbinom{n}{j}_q (-1)^j a^j q^{\frac{j(j-1)}{2}}.
\end{equation}

Next, we recall two identities due to Euler, which state (see (2.2.5) and (2.2.6) in \cite{And98})
that
\begin{align}
\label{E:Euler1}
\frac{1}{(x;q)_\infty} & = \sum_{n \geq 0} \frac{x^n}{(q;q)_n}, \\
\label{E:Euler2}
(x;q)_\infty &= \sum_{n \geq 0} \frac{(-1)^n x^n q^{\frac{n(n-1)}{2}}}{(q;q)_n}.
\end{align}
We also require summation and transformation identities for hypergeometric $q$-series. Recall that if $r \geq 1$, the basic hypergeometric $q$-series is defined by
\begin{equation*}
{}_{r+1}\phi_r \left(\begin{array}{cccc} a_1, & a_2, & \dots, & a_{r+1} \\
& b_1, & \dots, & b_r \end{array}; q; t \right) :=
\sum_{n \geq 0} \frac{(a_1, \dots, a_{r+1}; q)_n \; t^n}{(q, b_1, \dots, b_r; q)_n}.
\end{equation*}
Finite summations are obtained if one of the $a_i$ is set to be a negative power of $q$, and the following such evaluations are useful.
\begin{lemma}
\label{L:2phi1}
For $n \in\N$,
\begin{align}
\label{E:qChuVan}
{}_2 \phi_1 \left(\begin{array}{cc} a, & q^{-n} \\ & c \end{array}; q;  \frac{c q^n}{a} \right)
& =\frac{\left(\frac{c}{a}\right)_n}{\left(c\right)_n}, \\
\label{E:qChu2}
{}_2 \phi_1 \left(\begin{array}{cc} a, & q^{-n} \\ & c \end{array}; q;  \frac{c q^{n-1}}{a} \right)
&=\frac{\left(\frac{c}{a}\right)_{n-1}}{\left(c\right)_{n-1}}
- \frac{\frac{c}{aq}(1-a)\left(\frac{c}{a}\right)_{n-1}}{(c)_n}.
\end{align}
\end{lemma}
\begin{remark*}
Equation \eqref{E:qChuVan} is known as the {\it $q$-Chu-Vandermonde identity} (see (II.6) in \cite{GR90}).
\end{remark*}
\begin{proof}
Both formulas are related to specializations of Heine's second transformation for $_2\phi_1.$ Specifically, the following identity is found on page 38 of \cite{And98}:
\begin{equation*}
{}_2 \phi_1 \left(\begin{array}{cc} a, & b \\ & c \end{array}; q; t \right)
= \frac{\left(\frac{c}{b}\right)_\infty (bt)_\infty}{(c)_\infty (t)_\infty}
{}_2 \phi_1 \left(\begin{array}{cc} \frac{abt}{c}, & b \\ & bt \end{array}; q;  \frac{c}{b} \right).
\end{equation*}
Setting $b = q^{-n}, t = \frac{cq^{n-i}}{a}$ gives the specialization
\begin{equation*}
{}_2 \phi_1 \left(\begin{array}{cc} a, & q^{-n} \\ & c \end{array}; q;  \frac{c q^{n-i}}{a} \right)
= \frac{\left(\frac{cq^{-i}}{a}\right)_n}{(c)_n} \,
{}_2 \phi_1 \left(\begin{array}{cc} q^{-i}, & q^{-n} \\ & \frac{cq^{-i}}{a} \end{array}; q; cq^n \right).
\end{equation*}
Equation \eqref{E:qChuVan} follows immediately from the case $i=0$, and \eqref{E:qChu2} follows from the case $i=1$ after regrouping and simplifying.
\end{proof}

We close this section with two recurrences satisfied by our family of double series (recall \eqref{E:definer}).
\begin{lemma}
\label{L:recurrences}
We have
\begin{align}
\label{E:rec1}
r(\ell, u, v, w; n) -r(\ell, u, v, w; n-1) &= q^n r(\ell+1, u, v, w; n), \\
\label{E:rec2}
r(v, u, v, w; n)-r(0, u, v, w; n)&=-q^{w-uv} r(0, u, v, w; n-u).
\end{align}
\end{lemma}
\noindent Both of these are straightforward to verify, although we note that the second recurrence requires the summation on the right-hand side to be re-indexed, sending $j \mapsto j-1$.

\section{Combinatorial arguments}
\label{S:Comb}

In this section we present combinatorial proofs of our double series identities. A similar argument was used in Section 3.2 of \cite{BLM13} to reprove the double series representation for partitions without $k$-sequences that was originally shown in \cite{And05}. The general philosophy of this approach uses gap conditions and modular diagrams to decompose partitions into independent combinatorial components, each of whose generating functions are simple summations.  These pieces are then put back together by combining the generating functions in the appropriate manner. We also note that a more intricate version of this procedure was used in \cite{AB05} to prove \eqref{E:AB1} and \eqref{E:AB2}.

If $\lambda$ is a partition and $d \geq 1$, then the {\it $d$-modular diagram} of $\lambda$ is a row diagram in which in the $i$-th row consists of a sequence $dd\cdots dr$ with $\lfloor \frac{\lambda_i}{d}\rfloor$ $d$s and $r = \lambda_i \mod{d}$. For example, the $3$-modular diagram of $30 + 26 + 23 + 18 + 12 + 8 + 4 + 1$ is
\begin{align*}
& 3333333333 \\
& 333333332 \\
& 33333332 \\
& 333333 \\
& 3333 \\
& 332 \\
& 31 \\
& 1
\end{align*}

\subsection{Proof of Theorem \ref{T:fSdouble}}
\label{S:Comb:S}
In the case of Schur partitions, we use $3$-modular diagrams. Let $\S^n$ denote the set of Schur partitions with $n$ parts.
In general, since each row in $\lambda$ is at least $3$ larger than the previous, each row in the $3$-modular diagram has a distinct length, so we can remove an upper-left triangle of size $\ell(\lambda) - 1$.  Denote the resulting partition by $\lambda'$, where
\begin{equation*}
\lambda'_i := \lambda_i - 3(\ell(\lambda) - i).
\end{equation*}
Clearly we have the relations
\begin{align}
\label{E:lambda'=lambda}
\ell(\lambda') &= \ell(\lambda), \\
|\lambda'| & = |\lambda| - 3 \frac{\ell(\lambda)(\ell(\lambda) - 1)}{2}. \notag
\end{align}

The conditions on Schur partitions show that $\lambda'$ is a partition in which the multiples of $3$ must be distinct, but the parts that are $1, 2 \pmod{3}$ are arbitrary. Denote the set of such partitions by $\calT$.  It is clear that the map $\lambda \mapsto \lambda'$ is a bijection from $\S$ to $\calT$ that preserves the number of parts, and we therefore have the generating function
\begin{equation}
\label{E:gamma'prod}
\sum_{\lambda \in \S} x^{\ell(\lambda')} q^{|\lambda'|}
= \sum_{\sigma \in \calT} x^{\ell(\sigma)} q^{|\sigma|}
= \frac{(-xq^3; q^3)_\infty}{(xq, xq^2; q^3)_\infty} = \frac{(x^2 q^6; q^6)_\infty}{(xq; q)_\infty}.
\end{equation}
Using \eqref{E:Euler1} and \eqref{E:Euler2} to expand the final products, this is equivalent to
\begin{equation}
\label{E:lambda'sum}
\sum_{\lambda \in \S} x^{\ell(\lambda')} q^{|\lambda'|} = \sum_{k, m \geq 0} \frac{(-1)^k x^{m + 2k} q^{m + \frac{6k(k+1)}{2}}}{(q;q)_m (q^6; q^6)_k}.
\end{equation}

Finally, we use \eqref{E:lambda'=lambda} to add back in the missing triangles to \eqref{E:lambda'sum}.  In particular, \eqref{E:lambda'sum} allows us to pick off the $\lambda$ with exactly $n$ parts, implying that
\begin{equation*}
\sum_{\lambda \in \S^n} q^{|\lambda'|} =
\sum_{m + 2k = n} \frac{(-1)^k q^{m + \frac{6k(k+1)}{2}}}{(q;q)_m (q^6; q^6)_k}.
\end{equation*}
Adding back in the 3-modular triangle of size $n-1$ gives
\begin{equation*}
\sum_{\lambda \in \S^n} q^{|\lambda|} = q^{\frac{3n(n-1)}{2}} \sum_{\lambda \in \S^n} q^{|\lambda'|}
= \sum_{m + 2k = n} \frac{(-1)^k q^{m + \frac{6k(k+1)}{2}} q^{ \frac{3(m+2k)(m+2k-1)}{2}}}{(q;q)_m (q^6; q^6)_k},
\end{equation*}
because $n = m + 2k$ in the sum.
The proof is complete upon summing over all $n$, as
\begin{equation}
\label{E:fSlambda}
f(x) = \sum_{\lambda \in \S} x^{\ell(\lambda')} q^{|\lambda'| + 3\frac{\ell(\lambda')(\ell(\lambda') - 1)}{2}}
=
\sum_{m, k \geq 0} \frac{(-1)^k x^{m + 2k} q^{m + \frac{6k(k+1)}{2}}q^{3 \frac{(m+2k)(m+2k-1)}{2}} }{(q;q)_m (q^6; q^6)_k}.
\end{equation}
This simplifies to give the theorem statement.

\begin{remark*}
For example, if $\lambda = 30 + 26 + 23 + 18 + 12 + 8 + 4 + 1 \in \S$, then $\ell(\lambda) = 7$, and $\lambda' = 9 + 8 + 8 + 6 + 3 + 2 + 1 + 1$, which does indeed have the properties described above.
\end{remark*}

\subsection{$2$-modular diagrams and the proof of Theorem \ref{T:GG}}
\label{S:Comb:GG}

We now briefly sketch the proof of Theorem \ref{T:GG}, which follows a very similar argument, the main difference is that we now use $2$-modular diagrams. Beginning with \eqref{E:GGdouble}, we note that if $\lambda \in \GG$, then each row in $\lambda$ is at least $2$ larger than the previous. The rows in the $2$-modular diagram therefore have distinct lengths, so we can remove an upper-left triangle of size $\ell(\lambda) - 1$.

The resulting partition $\lambda'$ has distinct even parts, but odd parts may be repeated, so the generating function is
\begin{align}
\label{E:2gamma'}
\sum_{\lambda \in \GG} x^{\ell(\lambda')} q^{|\lambda'|}
& = \frac{(-xq^2; q^2)_\infty}{(xq; q^2)_\infty} = \frac{(x^2 q^4; q^4)_\infty}{(xq; q)_\infty}
= \sum_{k, m \geq 0} \frac{(-1)^k x^{m + 2k} q^{m + \frac{4k(k+1)}{2}}}{(q;q)_m (q^4; q^4)_k}.
\end{align}
The final equality follows from \eqref{E:Euler1} and \eqref{E:Euler2}. Adding back in the $2$-modular triangle gives
\begin{equation*}
\label{E:gxlambda}
f_{\GG}(x;q) = \sum_{\lambda \in \GG} x^{\ell\left(\lambda'\right)} q^{\left|\lambda'\right| + 2\frac{\ell\left(\lambda'\right)\left(\ell\left(\lambda'\right) - 1\right)}{2}}
=
\sum_{k, m \geq 0} \frac{(-1)^k x^{m + 2k} q^{m + \frac{4k(k+1)}{2} + (m+2k)(m+2k-1)} }{(q;q)_m \left(q^4; q^4\right)_k},
\end{equation*}
which implies the first part of the theorem statement.

The combinatorial proof of \eqref{E:Gdouble} is analogous, except now $\lambda'$ has arbitrary even parts and distinct odd parts, so the third expression in \eqref{E:2gamma'} is replaced by $\frac{\left(x^2 q^2; q^4\right)_\infty}{\left(xq; q\right)_\infty}.$ This implies the second part of the theorem.

\section{Analytic arguments I: $q$-differences and recurrences}
\label{S:AnalyticI}

In this section we give our first set of analytic proofs for the double series identities in Theorems \ref{T:fSdouble} -- \ref{T:Rindentities}. We use the theory of $q$-difference equations and recurrence relations.

\subsection{Proof of Theorem \ref{T:fSdouble}}
\label{S:AnalyticI:S}
We begin with a known $q$-difference equation for Schur partitions. Equation (2.11) of \cite{And68} states that
\begin{equation}
\label{E:fSdiff}
f_\S(x) = \left(1 + xq + xq^2 \right) f_\S\left(xq^3\right) + xq^3\left(1-xq^3\right) f_\S\left(xq^6\right);
\end{equation}
the combinatorial properties of this equation are also studied in Proposition 2.1 of \cite{BM13}. If we expand the generating function as a series in $x$, writing
\begin{equation*}
f_\S(x) = \sum_{n \geq 0} x^n U_\S(n),
\end{equation*}
then \eqref{E:fSdiff} is equivalent to the recurrence
\begin{equation}
\label{E:Rrec}
\left(1-q^{3n}\right) U_\S(n) = \left(q^{3n-2} + q^{3n-1} + q^{6n-3}\right) U_\S(n-1)
- q^{6n-6} U_\S(n-2),
\end{equation}
with initial conditions $U_\S(0) = 1$ and $U_\S(n) = 0$ for all $n < 0$ (note that we are allowed to consider negative-indexed terms when working with hypergeometric $q$-series due to \eqref{E:q_n<0}).
If we replace $n$ by $n-1$, multiply the resulting recurrence by $q^{3n}$ and then subtract it from the above, we obtain the new recurrence
\begin{align}
\label{E:Rreclong}
\left(1 - q^{3n}\right) U_\S(n) & = q^{3n-2} \left(1 + q + q^2\right) U_\S(n-1) \\
&
\quad - q^{6n-6} \left(1 + q + q^2 + q^{3n-3}\right) U_\S(n-2) + q^{9n-12} U_\S(n-3), \notag
\end{align}
with the same initial conditions as before. Either \eqref{E:Rrec} or \eqref{E:Rreclong} are sufficient to uniquely determine the $U_\S(n)$, but it is the lengthened recurrence that we use in order show that the double series representation for $f_\S(x)$ is correct.

This lengthened recurrence \eqref{E:Rreclong} is also equivalent to the $q$-difference equation
\begin{align}
\label{E:fSdifflong}
f_\S(x) = & \left(1 + xq + xq^2 + xq^3\right) f_\S\left(xq^3\right)
- \left(x^2 q^6 + x^2 q^7 + x^2 q^8\right) f_\S\left(xq^6\right) \\
& + \left(x^3 q^{15} - x^2 q^9\right) f_\S\left(xq^9\right). \nonumber
\end{align}
We note that this equation can also be proven directly from the definition of Schur partitions through combinatorial inclusion-exclusion, as in the proof of Proposition 2.1 of \cite{BM13}. In particular, the first term encodes the fact that the smallest part in a Schur partition may be 1, 2, or 3, with the next part larger than 3. The second term encodes the fact that the following pairs of successive parts are excluded, respectively: $(2,4), (3,4),$ and $(3,5).$  The pair (3,6) is also excluded, which is encoded by $-x^2 q^9 f(xq^9)$. Finally, we have now doubly subtracted the excluded triple $(3,5,7)$, which is added back in by $x^3 q^{15} f(xq^9)$.

For $\ell,n \in \Z$, define
\begin{equation*}
\label{E:rho}
\rho_\ell(n) := r(\ell,2,3,4;n) = \sum_{j \geq 0} \frac{(-1)^j q^{3j^2 + j - 2\ell j}}{(q;q)_{n-2j} (q^6; q^6)_j}.
\end{equation*}
It is clear that $\rho_\ell(0) = 1$ and $\rho_\ell(n) = 0$ for $n < 0$.

Recalling \eqref{E:fSlambda} and substituting $n = m + 2k$, we see that Theorem \ref{T:fSdouble} is equivalent to
\begin{equation}
\label{E:fS=rho}
f_\S(x) = \sum_{n \geq 0} x^n q^{\frac{n(3n-1)}{2}} \rho_0(n).
\end{equation}
Observe that Lemma \ref{L:recurrences} implies the following system of
recurrences:
\begin{align}
\label{E:rho0}
\rho_0(n) - \rho_0(n-1) & = q^n \rho_1(n), \\
\label{E:rho01}
\rho_0(n-1) - \rho_0(n-2) & = q^{n-1} \rho_1(n-1), \\
\label{E:rho02}
\rho_0(n-2) - \rho_0(n-3) & = q^{n-2} \rho_1(n-2), \\
\label{E:rho1}
\rho_1(n) - \rho_1(n-1) & = q^n \rho_2(n), \\
\label{E:rho11}
\rho_1(n-1) - \rho_1(n-2) & = q^{n-1} \rho_2(n-1), \\
\label{E:rho2}
\rho_2(n) - \rho_2(n-1) & = q^n \rho_3(n), \\
\label{E:rho3}
\rho_3(n) - \rho_0(n) & = -q^{-2} \rho_0(n-2).
\end{align}
This is a non-degenerate linear system of $7$ equations in the $7$ variables $\rho_0(n)$, $\rho_1(n)$, $\rho_1(n-1)$, $\rho_1(n-2)$, $\rho_2(n)$, $\rho_2(n-1)$,  and $\rho_3(n)$, so there is a unique solution for $\rho_0(n)$ in terms of the ``constants'' $\rho_0(n-1), \rho_0(n-2),$ and $\rho_0(n-3).$

Indeed, it is not difficult to work out through direct substitution, as we have
\begin{align*}
\rho_0(n) & = \rho_0(n-1) + q^n \rho_1(n) \\
& = \rho_0(n-1) + q^n \left(\rho_1(n-1) + q^n \rho_2(n)\right) \\
& = \rho_0(n-1) + q\left(\rho_0(n-1) - \rho_0(n-2)\right)
+ q^{2n} \left(\rho_2(n-1) + q^n \rho_3(n)\right) \\
& = (1 + q) \rho_0(n-1) - q \rho_0(n-2)
+ q^{n+1} \left(\rho_1(n-1) - \rho_1(n-2)\right) \\
& \qquad \qquad + q^{3n} \left(\rho_0(n) - q^{-2}\rho_0(n-2)\right) \notag \\
& = q^{3n} \rho_0(n) + (1+q) \rho_0(n-1) - \left(q + q^{3n-2}\right) \rho_0(n-2) \\
& \qquad + q^2 \left(\rho_0(n-1) - \rho_0(n-2)\right)
- q^3 \left(\rho_0(n-2) - \rho_0(n-3)\right).
\end{align*}
Regrouping, we conclude the recurrence
\begin{align}
\label{E:rhoSrec}
\left(1 - q^{3n}\right) & \rho_0(n) \\
& = \left(1 + q + q^2\right) \rho_0(n-1)
- q\left(1 + q + q^2 + q^{3n-3}\right) \rho_0(n-2) + q^3 \rho_0(n-3). \notag
\end{align}

Recalling \eqref{E:fS=rho}, we set
\begin{equation*}
u_\S(n) := q^{\frac{n(3n-1)}{2}} \rho_0(n),
\end{equation*}
so that we want to show that $f(x) = \sum_n x^n u_\S(n).$ Using \eqref{E:rhoSrec}, we find an equivalent recurrence for the $u_\S(n)$, namely
\begin{align*}
\left(1 - q^{3n}\right) u_\S(n) = & q^{3n-2} \left(1 + q + q^2\right) u_\S(n-1)  \\
& \quad - q^{6n-6} \left(1 + q + q^2 + q^{3n-3}\right) u_\S(n-2)
+ q^{9n-12} u_\S(n-3).
\end{align*}
This is equivalent to \eqref{E:Rreclong}, and since the initial conditions of $u_\S(n)$ and $U_\S(n)$ also coincide, we have therefore verified that $u_\S(n) = U_\S(n)$ for all $n \geq 0$.  This completes the proof of Theorem \ref{T:fSdouble}.

\begin{remark*}
It is striking that the analytic proof only recovers the longer recurrence \eqref{E:Rreclong}, and not the shorter one \eqref{E:Rrec} that arises directly from \eqref{E:fSdiff}.
\end{remark*}

\subsection{Proof of Theorem \ref{T:GG}}
\label{S:AnalyticI:G}

We next turn to the G\"ollnitz-Gordon identities and briefly outline an analytic proof of Theorem \ref{T:GG} using $q$-difference equations. As the proof is very similar to Section \ref{S:AnalyticI:S}, we suppress the calculations that are analogous to those seen between \eqref{E:fSdiff} -- \eqref{E:fSdifflong}. This is reflected in the modified notation that we use in this section, as we recall \eqref{E:R=GG} in setting
\begin{equation}
\label{E:grExpn}
\widetilde{f}_\GG(x) = \widetilde{f}_\GG(x;q) := \sum_{n \geq 0} x^n q^{n^2} r(0, 2, 2,2; n),
\end{equation}
and further write $\rho_\ell(n) := r(\ell, 2, 2, 2;n)$.

The recurrences in Lemma \ref{L:recurrences} then imply
\begin{align}
\label{E:GGrhorec}
\rho_0(n) &= \rho_0(n-1) + q^n \rho_1(n) \\
& = \rho_0(n-1) + q^n\left(\rho_1(n-1) + q^n\rho_2(n)\right) \notag \\
& = \rho_0(n-1) + q\left(\rho_0(n-1) - \rho_0(n-2)\right) +
q^{2n}\left(\rho_0(n) - q^{-2}\rho_0(n-2)\right), \notag
\end{align}
where the final parenthetical grouping follows from \eqref{E:rec2}, and all other steps use instances of \eqref{E:rec1}. Regrouping, this simplifies to
\begin{equation*}
\left(1 - q^{2n}\right) \rho_0(n) = (1 + q)\rho_0(n-1)
- \left(q + q^{2n-2}\right) \rho_0(n-2).
\end{equation*}
Letting $u_\GG(n) := q^{n^2} \rho_0(n)$, we then find the recurrence
\begin{equation*}
\left(1-q^{2n}\right) u_\GG(n) = q^{2n-1} (1+q) u_\GG(n-1) - q^{4n-3} \left(1 + q^{2n-3}\right) u_\GG(n-2),
\end{equation*}
which, by \eqref{E:grExpn}, is equivalent to the $q$-difference equation
\begin{equation}
\label{E:fGGdifflong}
\widetilde{f}_\GG(x) = \left(1 + xq + xq^2\right) \widetilde{f}_\GG\left(xq^2\right) - x^2q^5 \widetilde{f}_\GG\left(xq^4\right)
- x^2 q^6 \widetilde{f}_\GG\left(xq^6\right).
\end{equation}

As with \eqref{E:fSdiff} and \eqref{E:fSdifflong}, it is straightforward to verify that $f_\GG(x)$ satisfies a simple $q$-difference equation, where we write $f_\GG(x) := f_\GG(x;q)$:
\begin{equation}
\label{E:fGGdiff}
f_\GG(x) = (1+xq) f_\GG\left(xq^2\right) + xq^2 f_\GG\left(xq^4\right).
\end{equation}
For example, this can be verified by expanding \eqref{E:Gsum} as a series in $x$; the formula also appears as (3.2) in \cite{And67}. Furthermore, one can understand the recurrence combinatorially by noting that the first term on the right generates those G\"ollnitz-Gordon partitions whose smallest part is either 1 or at least $3$ (so the next part is at least $3$), while the second term generates those whose smallest part is $2$ (so the next part is at least $5$).

Then \eqref{E:fGGdifflong} is obtained by substituting $x \mapsto xq^2$ in \eqref{E:fGGdiff}, multiplying the entire equation by $xq^2$, and subtracting the result from \eqref{E:fGGdiff} itself. As in Section \ref{S:AnalyticI:S}, the fact that $f_\GG$ and $\widetilde{f}_\GG$ satisfy the same recurrence is enough to conclude \eqref{E:GGdouble}.

Next, recall \eqref{E:R=G} and set
\begin{equation*}
\widetilde{f}_{\G}(x) := \sum_{n\geq 0} x^n q^{n^2} r\left( 1,2,2,2;n\right).
\end{equation*}
Noting that $r(1,2,2,2;n) = r(0,2,2,0;n)$, we also set $\rho_j(n) := r(j,2,2,0;n)$. Proceeding as above, we obtain a formula that is identical to \eqref{E:GGrhorec} except that the final term is now $q^{-4}\rho_0(n-2)$. This translates to the functional equation
\begin{equation}
\label{E:fGdifflong}
\widetilde{f}_\G(x) = \left( 1+xq +xq^2\right) \widetilde{f}_\G\left( xq^2 \right) - x^2 q^5 \widetilde{f}_\G\left( xq^4\right) - x^2 q^4 \widetilde{f}_\G\left( x q^6\right).
\end{equation}

Similar to the above, one easily finds the recurrence
\begin{equation*}
f_\G(x) = \left( 1+xq^2\right) f_\G\left( xq^2\right) +x q f_\G\left( xq^4\right),
\end{equation*}
which, upon comparison with \eqref{E:fGdifflong}, concludes the proof of \eqref{E:Gdouble} as before.

\subsection{Proof of Theorem \ref{T:Rindentities}}
\label{S:AnalyticI:R}

We conclude this section by using $q$-difference equations to prove the identities in Theorem \ref{T:Rindentities}. We simultaneously prove \eqref{E:R=Euler1} and \eqref{E:R=Euler2}, by first setting $\rho_j (n):= r(j,1,2,1;n).$ Using Lemma \ref{L:recurrences}, we find that
\begin{equation}
\label{E:rhoEuler1}
\left( 1-q^{2n} \right) \rho_0 (n) = \left( 1+q-q^{2n-1}\right)\rho_0 (n-1) - q \rho_0 (n-2).
\end{equation}
Following the arguments in the previous subsections, \eqref{E:rhoEuler1} implies that the series $f_E(x) = f_E(x;q) := \sum_{n \geq 0} x^n q^{n^2} \rho_0(n)$ satisfies the $q$-difference equation
\begin{equation}
\label{E:fEdifflong}
f_E(x) - f_E \left( xq^2\right) = \left( 1+q\right) xq f_E\left( xq^2\right) - q^2 xf_E\left( xq^4\right)- x^2 q^5 f_E\left( xq^4\right).
\end{equation}
An short calculation shows that $(-xq; q^2)_\infty$ also satisfies \eqref{E:fEdifflong}, so $f_E(x) = (-xq; q^2)_\infty.$ Since by definition $f_E(q^{t-1};q) = R(2,t,1,1,2,2)$, we conclude \eqref{E:R=Euler1} and \eqref{E:R=Euler2}, respectively, by setting $x = 1$ and $q$, respectively.

\begin{remark}
We have not included a proof of \eqref{E:Rm=prod} in this section as the resulting recurrences are essentially trivial (cf. Section \ref{S:AnalyticII:R}).
\end{remark}

\section{Analytic arguments II: Summation identities}
\label{S:AnalyticII}

In this section we use well-known $q$-series summation formulas in order to give analytic proofs for many of the double-series identities in Theorems \ref{T:GG} and \ref{T:Rindentities}. Notably, we have not found any such proof for Theorem \ref{T:fSdouble}, which seems to be more novel.

\subsection{Proof of Theorem \ref{T:GG}}
\label{S:AnalyticII:GG}

Recalling \eqref{E:R=GG} and \eqref{E:R=G}, we begin by noting that for $u \in\N_0$,
\begin{equation}
\label{E:ru222}
r(u,2,2,2;n) = \sum_{j \geq 0} \frac{(-1)^j q^{2j^2-2uj}}{(q)_{n-2j}\left(q^4;q^4\right)_j}.
\end{equation}
We rewrite this expression using \eqref{E:qnegn}, which implies that
\begin{equation*}
\frac{1}{(q)_{n-2j}}
=\frac{q^{2nj-j(2j-1)}\left(q^{-n}\right)_{2j}}{(q)_n}=
\frac{q^{2nj-j(2j-1)}\left(q^{-n}, q^{1-n}; q^2\right)_j}{(q)_n}.
\end{equation*}
Plugging in to \eqref{E:ru222}, we therefore have the hypergeometric expression
\begin{align}
\label{E:ru222=2phi1}
r(u, 2, 2, 2; n) &= \frac1{(q)_n}\sum_{j\geq 0}\frac{\left(q^{-n}, q^{1-n}; q^2\right)_j}{\left(q^2, -q^2; q^2\right)_j}\left(-q^{1+2n-2u}\right)^j \\
& = \frac{1}{(q)_n} {}_2 \phi_1 \left(\begin{array}{cc} q^{-n}, & q^{1-n} \\ & -q^2 \end{array}; q^2;  -q^{1+2n-2u} \right). \notag
\end{align}

For the G\"ollnitz-Gordon partitions, we set $u=0$, and then further distinguish cases based on the parity of $n$. If $n$ is even, then we use \eqref{E:qChuVan} with $n\mapsto\frac{n}{2}$, $q\mapsto q^2$, $a=q^{1-n}$, and $c = -q^2$ to get
\begin{equation}
\label{E:r0222neven}
r(0,2,2,2;n) = \frac1{(q)_n} \frac{\left(-q^{1+n}; q^2\right)_{\frac{n}{2}}}{\left(-q^2; q^2\right)_{\frac{n}{2}}}=\frac{\left(-q; q^2\right)_n}{(q)_n \left(-q^2; q^2\right)_{\frac{n}{2}}\left(-q; q^2\right)_{\frac{n}{2}}}
=\frac{\left(-q; q^2\right)_n}{\left(q^2; q^2\right)_n}.
\end{equation}
If $n$ is odd, then we instead use \eqref{E:qChuVan} with $n\mapsto \frac{n-1}{2}$, $q\mapsto q^2$, $a=q^{-n}$, and $c=-q^2$  to get
\begin{equation}
\label{E:r0222nodd}
r(0, 2, 2, 2; n)=\frac1{(q)_n}\frac{\left(-q^{2+n}; q^2\right)_{\frac{n-1}{2}}}{\left(-q^2; q^2\right)_{\frac{n-1}{2}}}=\frac{\left(-q; q^2\right)_n}{(q)_n\left(-q^2; q^2\right)_{\frac{n-1}{2}}\left(-q; q^2\right)_{\frac{n+1}{2}}}
=\frac{\left(-q; q^2\right)_n}{\left(q^2; q^2\right)_n}.
\end{equation}
As these evaluations simplify to a uniform expression, we find that
\begin{equation*}
f_\GG(x;q) =
\sum_{n\geq 0}\frac{x^n q^{n^2}\left(-q; q^2\right)_n}{\left(q^2; q^2\right)_n} =
\sum_{n\geq 0} x^n q^{n^2} r(0,2,2,2;n),
\end{equation*}
which proves \eqref{E:GGdouble}.

For the G\"ollnitz partitions we set $u=1$ in \eqref{E:ru222=2phi1}, and again separate cases by the parity of $n$. If $n$ is even, then we apply \eqref{E:qChu2} with $n \mapsto \frac{n}{2}$, $q \mapsto q^2$, $a = q^{1-n}$, and $c=-q^2$ to obtain
\begin{align*}
{}_2 \phi_1 \left(\begin{array}{cc} q^{-n}, & q^{1-n} \\ & -q^2 \end{array}; q^2;  -q^{2n-1} \right)
& = \frac{\left(-q^{n+1};q^2\right)_{\frac{n}{2}-1}}{\left(-q^2; q^2\right)_{\frac{n}{2}-1}}
+ \frac{q^{n-1} \left(1 - q^{1-n}\right) \left(-q^{n+1}; q^2\right)_{\frac{n}{2} - 1}}{\left(-q^2; q^2\right)_{\frac{n}{2}}} \\
& = \frac{q^n \left(1 + q^{-1}\right) \left(-q^{n+1}; q^2\right)_{\frac{n}{2} - 1}}{\left(-q^2; q^2\right)_{\frac{n}{2}}}.
\end{align*}
This further simplifies as in \eqref{E:r0222neven}, giving the overall expression
\begin{equation*}
r(1,2,2,2;n) = \frac{q^n \left(-q^{-1}; q^2\right)_n}{\left(q^2; q^2\right)_n}.
\end{equation*}

If $n$ is odd, then we apply \eqref{E:qChu2} with $n \mapsto \frac{n-1}{2}$, $q \mapsto q^2$, $a = q^{-n}$, and $c=-q^2$, obtaining
\begin{align*}
{}_2 \phi_1 \left(\begin{array}{cc} q^{-n}, & q^{1-n} \\ & -q^2 \end{array}; q^2;  -q^{2n-1} \right)
& = \frac{\left(-q^{n+2};q^2\right)_{\frac{n-3}{2}}}{\left(-q^2; q^2\right)_{\frac{n-3}{2}}}
+ \frac{q^{n} \left(1 - q^{-n}\right) \left(-q^{n+2}; q^2\right)_{\frac{n-3}{2}}}{\left(-q^2; q^2\right)_{\frac{n-1}{2}}} \\
& = \frac{q^n \left(1 + q^{-1}\right) \left(-q^{n+2}; q^2\right)_{\frac{n-3}{2}}}{\left(-q^2; q^2\right)_{\frac{n-1}{2}}}.
\end{align*}
Again the overall expression for $r(1,2,2,2;n)$ simplifies further (cf. \eqref{E:r0222nodd}), giving a uniform expression in $n$. This completes the proof of \eqref{E:Gdouble}.

\subsection{Proof of Theorem \ref{T:Rindentities}}
\label{S:AnalyticII:R}

We begin with \eqref{E:R=Euler1} and \eqref{E:R=Euler2}. By the definition of $r(\ell, u, v, w; n)$ and \eqref{E:qnegn} and \eqref{E:qtlim}, we find the hypergeometric limit
\begin{align*}
r(1,1,2,2;n) & =  \frac{1}{(q)_n}\lim_{t \rightarrow 0} \sum_{j\geq0} \frac{\left(q^{-n}\right)_j \left(-\frac{q}{t}\right)_j}{(q)_j (-q)_j}t^j q^{nj}
=\frac{1}{(q)_n} \lim_{t\rightarrow 0}
{}_2\phi_{1}\left(\begin{array}{cc}q^{-n} ,& -\frac{q}{t} \\ & -q \end{array}; q; tq^n \right).
\end{align*}
Using \eqref{E:qChuVan} with $a = -\frac{q}{t}$ and $c = -q$, this evaluates to
\begin{equation*}
\frac{1}{(q)_{n}} \lim_{t\rightarrow 0}\frac{(t)_{n}}{(-q)_n}=\frac{1}{(q^2;q^2)_n}.
\end{equation*}
Euler's identity \eqref{E:Euler2} now implies the two desired formulas,
\begin{align*}
 R(2,1,1,1,2,2) = \sum_{n\geq 0} \frac{q^{n^2}}{(q^2;q^2)_n}  &=\left(-q;q^2\right)_{\infty}, \\
 R(2,2,1,1,2,2) = \sum_{n\geq 0} \frac{q^{n^2+n}}{(q^2;q^2)_n}&=\left(-q^2;q^2\right)_{\infty}.
\end{align*}

Finally, for \eqref{E:Rm=prod}, we observe that if we set $a = q$ in \eqref{E:qbinom} and divide by $(q)_n$, this immediately implies that $r(1,1,1,2;n) = 1$. We therefore have the formula
\begin{equation*}
R(m,m,1,1,1,2) = \sum_{n \geq 0} q^{\frac{mn(n+1)}{2}} = \frac{\left(q^{2m}; q^{2m}\right)_\infty}{\left(q^{m}; q^{2m}\right)_\infty},
\end{equation*}
where the final equality follows from Gauss' Triangular Number identity
((2.2.13) in \cite{And98}).

\section{Further questions}
\label{S:Conc}

We close with a small collection of additional questions that naturally arise from our study of double series representations for combinatorial partition identities.

\begin{enumerate}

\item
Schur's Theorem
also extends to partitions where the parts differ by at least $d$ for $d \geq 3$. Indeed, there is a two-parameter family of such identities indexed by $(d,r)$, and the analytic properties of these generating functions were recently described in work of the second and third author \cite{BM13}.

However, the techniques of this paper do not give a double series
for Schur's partitions if $d \geq 4$. For example, the combinatorial proof of Theorem \ref{T:fSdouble} in Section \ref{S:Comb:S} yields a hypergeometric series with a triple summation in general, as the final simplifications in \eqref{E:gamma'prod} only work for $d=3$. Furthermore, it is known that the generating function for Schur's general $(d,r)$-family also satisfy $q$-difference equations similar to \eqref{E:fSdiff} (see formula (2.1) of \cite{AG93} or Proposition 2.1 of \cite{BM13}). As far as we have been able to determine, these $q$-difference equations are not compatible with the recurrences from Lemma \ref{L:recurrences}, which is another indication of the obstacles to proving double series representations.

\item
One also sees $3$-modular diagrams in Pak's elegant bijective proof of Schur's identities in 4.5.1 and 4.5.2 of \cite{Pak06}.  It would be interesting to determine if there is any relationship between his map and the $3$-modular decomposition used in the combinatorial proof of \ref{T:fSdouble} from Section \ref{S:Comb:S}.

\end{enumerate}

\end{document}